\documentclass[11pt]{article}
\usepackage{amsmath,amssymb,mathrsfs}
\usepackage[ansinew]{inputenc}
\usepackage{makeidx}
\makeindex
\advance\topmargin 0cm
\def\Q{\mbox{\small Q}}
\def\R{\mathbb{R}}
\def\C{\mathbb{C}}
\def\Bn{\mathbb{B}^n}

\def\i{\imath}
\def\fl{\rightarrow}

\def\B{{\cal B}}

\def\dst{\displaystyle}

\def\th {\tanh}
\long\def\cor #1#2 {\par{\bf Corollary #1. }{\it #2}\par}
\long\def\lemme #1#2 {\par\noindent{\bf Lemma #1. }{\it #2}\par}
\long\def\prop #1#2 {\par\noindent{\bf Proposition #1. }{\it #2}\par}
\long\def\Th #1#2{\par\noindent{\bf #1. }{\it #2}\par}

\long\def\Lemme#1{\par\noindent{\bf Lemma. }{\it #1}\par}

\def\dem{\par\noindent{{\it Proof\,}: }}
\advance \parskip 0pt

\begin{document}
{\centering {\large\bf  Lie Ball as Tangent Space to Poincaré Ball}

\small by  Roger \sc Tchangang Tambekou \par}\vskip 12pt

{\abstract{We equip the whole tangent space $TM$ to
a hyperbolic manifold $M$ (of constant sectional curvature -1) with
a natural metric in an intrinsic way, so that the isometries of
$M$ extend to isometries of $TM$ by holomorphic continuation. The image to the tangent space to a geodesic is equivalent to a hyperbolic disk. 

In the case of hyperbolic space, we exhibit 
an equivariant diffeomorphism  between $TM$ and the fourth symmetric
complex domain of E. Cartan, also known as the Lie ball.
The closure of the Lie ball appears as a horospheric compactification of the tangent bundle to hyperbolic
space, and its Bergmann metric gives an intrinsic natural k\"ahler metric
on the tangent space $TM$.

 The equivariant map has a simple geometric interpretation.

}}
\medskip

We propose hereafter another complexification of the hyperbolic space, at least as 'natural' as the one given by the Akhiezer-Gindikin domain. The result is still the Lie ball, but here, not only the complex structure is global, but the leaves of the riemannian foliation are complete hyperbolic disks instead of having partial flat structure. This leads to an equivariant compactification of hyperbolic space.

We define a diffeomorphism between the Lie ball and the tangent space to the Poincaré ball, and pull back the Kähler symmetric structure of the Lie ball to the tangent space to hyperbolic space. In the sequel, $n$ denotes an integer greater or equal to 2.

\Th{Main result}{The action of the hyperbolic group extends to the closure $\overline{\B^n}$ of the Lie ball as isometries, by holomorphic continuation.

There exists a diffeomorphism $\theta$ between $T\Bn$, the tangent space to the hyperbolic space, and  the Lie ball $\B^n$, with the following properties:
\begin{enumerate}
  \item $\theta$ is equivariant for the action of the hyperbolic group.
  \item $\overline{\B^n}$ is an equivariant compactification of the tangent space $T\Bn$ by codimension 2 horospheres.
  \item The tangent space to every geodesic of the hyperbolic space is totally geodesic. Its image by $\theta$ is isometric to a hyperbolic disk.
\item Any vector line of the tangent space at a point is a geodesic.
\end{enumerate}}

These past twenty years, serious efforts have been undertaken to equip the tangent space to a
hyperbolic manifold with a k\"ahler structure. Several extensions to
$TM$ of the metric on $M$ have been studied. The most famous are the
Sasaki and the Cheeeger-Gromoll metrics, which are too rigid in the
sense that $TM$ is a complex space form only if $M$ is flat. More
recently, Oproiu [14], then Abbassi [1] gave a family of natural
metrics constructed following a technique due to Dombrowski, with a slight modification of
the Sasaki construction. In the case of a hyperbolic
manifold, they can turn $TM$ into a k\"{a}hler locally symmetric
space. But the differential geometry of such a space seems hard to
understand.

More powerful approaches come from Lie Group Theory. Following
Akhiezer and Gindikin [2], Burns, Halverscheid and Hind [7] compute the
canonical complexification of riemannian symmetric spaces of non
compact type, by using the technique of adapted complex structures
brought independently by Guillemin-Stenzel [9] and Lempert-Szöke [13].  For
 hyperbolic space, they find the Lie ball.

To define the so-called canonical complexification of a
riemannian symmetric space $M$, Akhiezer and Gindikin define a map
$\psi$ from the tangent space to the complexified of $M$ (in the
sense of Lie group theory). The map $\psi$ is a diffeomorphism in a
neighborhood of the 0-section. For every real $r>0$, let
$\Omega_r=\{v\in TM, \|v\|<r\}$ be the Grauert tube of radius $r$. The \emph{complexification} of $M$, also
called its complex crown is $\psi(\Omega_r)$ where $r$ is the
largest value for which $\psi$ is a diffeomorphism onto its
image.

In the case of the hyperbolic space, $r=\pi/2$, and the image of
$\Omega_r$ is the fourth symmetric domain of Elie Cartan, called  Lie ball. An interesting fact is that $\psi$  is equivariant w.r.t. the action of the hyperbolic group, i.e. the complex
analytic continuation of an isometry of $M$ acts the same way than
the tangent map on $\Omega_r\subset TM$.

It is obvious that such an equivariant action can be extended to the whole of $TM$, but at the expense
of loss of some essential properties.  The horizontal and the vertical $n$-planes of the second tangent space
are no more orthogonal. Every connection on $M$ determines such a
distribution.

So, unlike in the compact case, the technique of adapted structures does not yet lead to an interesting  complexification of the entire tangent space for non compact symmetric spaces.

Our method is entirely geometric and leads to an equivariant
diffeomorphism between $TM$ and the Lie ball, with a
natural metric for which the horizontal distribution and the
vertical one prescribed by the connection on $M$ are orthogonal.
Moreover, it gives
a  compactification of the tangent space to any complete hyperbolic manifold.

The tangent space to a geodesic of the hyperbolic space has the structure of a hyperbolic disk and is totally
geodesic.

We use a technical tool, the $T$-map,
introduced by Lelong [12] and
Aronsjahn to develop the theory of harmonicity cells, in view of the
study of singularities of complex extensions of harmonic
functions. It has some resemblances with the theory of linear cycles
in flag domains developed by Wolf and al. [16], but is by far
more elementary.

Part of the study carried here for
the unit ball generalizes to any domain of $\R^n$. We have
kept only geometric ideas, being as self-contained as possible. For
the reader interested in holomorphic continuation, see [3],[6],[12].

For hyperbolic geometry, chapter 2 of W. Thurston's book [15]
is the reference.

\noindent {\bf Acknowledgment. } I thank Noël Lohoue for his encouragements. He kept me informed of the Lie group specialists point of view on the subject. I also thank my colleagues from Yaounde I University Mathematics  Department. By discharging me from teachings, they gave me more time for research.

In the years 80, Avanissian renewed the interest of analysts in harmonicity
cells in his book dedicated to holomorphic continuation [3]. We are
deeply indebted to his teachings and our work benefited from numerous
remarks and critics from him.

We first associate to the space of $(n-2)$-spheres of the unit ball
two parametrizations, respectively by the  Lie ball $\B^n$ and the tangent space $T\Bn$.
We show that they are equivariant for the
induced action of the hyperbolic group. So, we can define a geometric transformation between the two
 spaces,
and deduce an intrinsic natural ([10],[11]) metric  on the whole of $TM$.

Let $\Bn$ be the Poincaré ball equipped with the hyperbolic metric
$\frac{ds^2}{(1-x^2)^2}$ at point $x$, where $ds^2$ is the euclidean
metric. The conformal models of hyperbolic geometry and Euclidean
geometry share a surprising common property: the spaces of
spheres are the same. But the euclidean center and radius are
different from hyperbolic ones, in general.

 Let $\mathcal{S}(\Bn)$ be the
space of $(n-2)$-spheres contained in $\Bn$.  We give two parametrizations of
$\mathcal{S}(\Bn)$, according to the geometry one chooses. To help develop
good mental pictures, think $n=3$.
$\R^n$ is oriented once for all.

An $(n-2)$-sphere $c$ contained in the unit ball has an euclidean
center $x$, an euclidean radius $r$, and is contained in a hyperplane. Let $y$ be a
vector orthogonal to that hyperplane, of length $r$.  We associate to $c$ the
complex number $z=x+\imath y$. To avoid confusion between $z$ and
$\overline{z}$, we proceed in the following way:
 let $b$ be a basis of the
tangent space to $c$ at point $p$, and $n$ the external normal at
$p$ to $c$ in the hyperplane containing it. The {\it radius vector} $y$ is chosen so that the basis $(y,n,b)$ of
$\R^n$ is direct.
This \textit{induced orientation}
 turns our parametrization into a bijective map between
$\C^n$ and the space $\mathcal{S}(\Bn)$ of (oriented) $(n-2)$-spheres of $\R^n$, for $n \geq 3$.
That map is called by Lelong [12] the $T$-map. (Lelong and the analysts did not need 
an orientation on the spheres.)
For $n=2$, one chooses the following
arbitrary order on the 0-sphere:\\
$T(z):=(z_1+\imath z_2, \
\overline{z_1}+\imath \overline{z_2})$ with $\mathbb{R}^2\sim \mathbb{C}$, as usual.

Let's now explain the (\emph{unusual}) action of maps on $\mathcal{S}(\Bn)$: in the sequel,
we are interested only in maps
preserving the $(n-2)$-spheres. Such maps are known as Moebius transformations.
 The image of a $(n-2)$ sphere $c$ is well defined setwise. Let $p\in c$. One chooses the radius
 vector $Y$ of $f(c)$
so that $Y\!\cdot df(p)(y)>0$. If $f$ preserves orientation, this
definition coincides with
the usual action of a differentiable transformation on an oriented submanifold.
If $f$ reverses orientation, there is a difference with the usual case. For instance, if we take
for $f$ the reflection through a coordinate hyperplane, and for $c$ an oriented $(n-2)$-sphere
 in the same hyperplane, then $f$ restricted to $c$ is the identity, but the
action of $f$ reverses the orientation of $c$.

We consider the quadratic form $\Q(z)=\sum_{1}^{n}z_i^2$ defined on $\C^n$, which extends the
euclidean real norm. It is obvious that
$$\forall a\in \R^n, z\in \C^n, \quad a\in T(z)\Leftrightarrow \Q(z-a)=0.$$
So, $a\in T(z)\Leftrightarrow z\in \Gamma(a)$, where $\Gamma(a)$ is the
$\Q$-isotropic cone with vertex $a$.

Let's compute the necessary and sufficient condition for a complex point $z=x+\i y$ in
$\C^n$ to correspond to a $(n-2)$-sphere contained in $\Bn$. If $y\neq 0$, any plane
$\varpi$ containing the origin $o$, $x$ and $y$,  is a symmetry plane
for the Poincaré ball and the $(n-2)$-sphere $T(z)$, which it
intersects in two points: the nearest to the origin, and $p$, the
furthest. Let $\rho$ be the rotation with angle $\pi/2$ in the plane
$\varpi$, of center $x$, which sends $x+y$ to $p$. Then $T(z)\subset
\Bn \Leftrightarrow\|p\|=\|x+\rho (y)\|<1$. Now,
\begin{eqnarray*}
  \|x+\rho(y)\|^2 &=& \parallel x\|^2+\|y\|^2 +
2\|x\|\|y\||\sin(x,y)| \\
   &=& \|x\|^2+\|y\|^2 +
2\sqrt{\|x\|^2\|y\|^2-{(x\cdot y)^2}}
\end{eqnarray*}

Therefore, the  condition  is:
$$\|z\|^2+\sqrt{\|z\|^4-\|\Q(z)\|^2}<1,$$
which is the equation of  the fourth symmetric domain of Elie Cartan,
 also called  Lie ball of complex dimension $n$, denoted here by ${\cal B}^n$. 
 
 The Lie ball $\B^n$ is the smallest convex domain stable by multiplication by unit complex numbers and containing the real unit ball. 
 It is included in the unit complex ball. Its
group of isometries is isomorphic to $SO(n,2)/\mathbb{Z}_2$ [8], which acts by holomorphic transformations. The
Bergmann metric [5] of the Lie ball is the only riemannian metric (up to a
constant positive factor) invariant under that action. 
It induces the hyperbolic metric on the Poincaré ball.

The stabilizer of the origin is isomorphic to the product 
$O(n)\times SO(2)$, where an isometry of $O(n)$ acts on the Lie
ball by holomorphic extension, and the action of $SO(2)=\{e^{\imath
t},\ t\in \mathbb{R}\}$ is mere multiplication by 
complex numbers.

The other parametrization of the space of $(n-2)$-spheres of $\Bn$
is obtained by considering the similar construction in hyperbolic
geometry. An $(n-2)$-sphere $c$ contained in the unit ball has a
hyperbolic center $x$, a hyperbolic radius $r$, and is contained
in a hyperbolic hyperplane, which is a sphere orthogonal to the
boundary of $\Bn$. To sphere $c$, we associate the
vector $v$ tangent to $\Bn$ at point $x$, of length $r$,
orthogonal to the hyperbolic hyperplane containing $c$. The
orientation of the space determines the choice between $v$ and
$-v$ exactly as before, in euclidean case.

This construction defines a bijective
map $S$ between the tangent space $T\Bn$ to the
Poincaré ball and the space of $(n-2)$-spheres $\mathcal{S}(\Bn)$
of the ball.
\Th{Lemma}{$S$ is an equivariant map with
respect to the induced actions of the hyperbolic group on the two
spaces.}
\dem Let $S(v)$ be the sphere associated to the tangent vector $v$ at point $x\in \Bn$, and $\varpi$
the hyperbolic hyperplane containing $S(v)$.
A hyperbolic motion $\rho$ transforms $S(v)$ into  a $(n-2)$-sphere centered at
$\rho(x)$ with same radius.
 Since $v$ is orthogonal to  $\varpi$ and
 $\rho$ is an isometry,  $d\rho( v)$ is orthogonal to the
 hyperplane
$\rho(\varpi)$ at $\rho(x)$. It follows that we have to choose the right image for $S(v)$
between $S(d\rho(v))$ and $S(-d\rho(v))$. So $\rho(S(v))=S(\varepsilon(\rho, v)d\rho(
v))$, with $\varepsilon(\rho, v)=\pm1$. One shows easily that $\varepsilon$ is a continuous map on the cartesian
product of the  space of hyperbolic isometries and the tangent space to the ball,
which has two connected components. For the identity map,
$\varepsilon(\rho, v)=1$, and  for the  reflection $\rho$
with respect to the hyperbolic hyperplane containing $S(v)$, one has $\varepsilon(\rho,
v)=-1$. We conclude that $\varepsilon(\rho, v)=\varepsilon(\rho)=1$ for
orientation preserving isometries and $-1$ for orientation
reversing ones.

\Th{Proposition} {the action of the
hyperbolic group on the Poincaré ball $\Bn$ extends to the Lie ball
$\B^n$ by holomorphic continuation.

The  bijection $T : \B^n \rightarrow \mathcal{S}(\Bn)$ is an
equivariant map between the Lie
ball and the space of oriented $(n-2)$-spheres of the Poincaré ball, all equipped
with the induced  action of the hyperbolic group.
}

Before proving the proposition, we first recall some useful properties of hyperbolic motions.

The hyperbolic group is generated by inversions with respect to
spheres orthogonal to the boundary sphere, called sphere
at infinity. An inversion is the hyperbolic analog of an euclidean
hyperplane symmetry. In fact, the orthogonal symmetries with respect to a hyperplane passing
through the origin are isometries in the two geometries and any
inversion is conjugated to them. So, they are orientation
reversing.
Let $a\in \R^n$.
The general equation of an inversion $\gamma$ with respect to
the sphere of center $a$ and radius $\alpha $ is

$$\gamma : \quad \begin{array}[t]{rcl}
\R^n-\{a\}&\fl &\R^n-\{a\}\\
x&\mapsto&\dst a+\alpha^2\frac{x-a}{\|x-a\|^2}\\
\end{array}$$
Such a map is involutive, i.e. equal to its inverse. Inversions are
better understood when defined on the Riemann sphere $S^n=\R^n \cup
\infty$. They are then diffeomorphisms. A Moebius map is a finite product
of inversions. It is a transformation of the sphere. The Moebius
group contains of course the euclidean isometries.

 A
differentiable map is said to be conformal on an open set of $\R^n$
if it preserves the non oriented euclidean angles. If a conformal
map reverses orientation, some call it anticonformal. 

A very important
characterization of Moebius maps is the following [4]:

\Th{Liouville's Theorem}{Let $\Omega$ be a domain in $\R^n$ and $f: \Omega \fl \R^n$ a
differentiable map.

\begin{enumerate}
\item For $n=2$, the conformal maps preserving the orientation are
the holomorphic functions with nonvanishing derivatives. Those that
reverse the orientation are the conjugates of those holomorphic
functions.
\item For $n\geq 3$, the only conformal maps are the
Moebius transformations.

 The conformal maps defined on
 the whole of $\R^n$ are 
compositions of affine homotheties  and euclidean isometries.
\end{enumerate}}

Let $\lambda$ be a Moebius map, acting on the Riemann sphere, such that $\lambda(\infty)\neq\infty$, and
$\gamma$  an inversion sending $a=\lambda^{-1} (\infty)$ to $\infty$. The
Moebius transformation $\lambda \circ \gamma^{-1}$ is conformal on
$\R^n$, so is equal to an affine map $\sigma$, according to
Liouville. So, $\lambda = \sigma \circ \gamma$. The holomorphic
continuation of $\lambda$ is then defined on $\C^n-\Gamma (a)$, where
$\Gamma (a)$ is the isotropic cone of $a$, i.e. the set of complex
points verifying $\Q(z-a)=0$. Thus,
\Th{Lemma}{Let $z\in \C^n$. If a Moebius map $\lambda$ is
defined on $T(z)$, then its holomorphic continuation is well defined
at $z$.}

It is well-known that the general form $\gamma$ of a hyperbolic motion of the Poincaré ball  is:
$$\gamma :\begin{array}[t]{rcl}
\B^n&\fl &\B^n\\
x&\mapsto&\dst \rho \frac{(\|a\|^2-1)x+(1+\|x\|^2)a-2(x\cdot a)a}{ \|x\|^2\|a\|^2-2x\cdot a+1}:=\rho \delta_a(x)
\end{array}$$
where $\rho$ is an euclidean isometry of the ball. The map
$\delta_a$ is an involution  interchanging $a\in
\Bn$ and the origin.

The denominator of $\gamma$ is $\|x\|^2\|a\|^2-2x\cdot a+1=
\|a\|^2\|x-a^*\|^2$, where $a^*=\frac{a}{\|a\|^2}$ is the image of
$a$ by the inversion with respect to the unit sphere, for
$a$ real. So, $\gamma$ has a holomorphic continuation in the Lie ball,
since for $z\in \mathcal{B}^n$, we have $T(z)\subset \Bn$, and
$a^*\not\in \Bn$.

To show that $T$ is equivariant with respect to the induced action
of the hyperbolic group, we will use the following striking
property of Moebius maps:
\Lemme {Let  $\gamma : \Bn \fl \Bn$
be a hyperbolic isometry
  and $\tilde \gamma $ its
holomorphic continuation, defined on the Lie ball $\B^n$.
For all $ z\in \B^n$,
\begin{itemize}
\item If $\gamma$ is orientation preserving, then $\gamma (T(z))
=T(\tilde \gamma (z))$
\item if  $\gamma$ is orientation
reversing, then $\gamma (T(z))=T(\overline{\tilde \gamma (z)})$
\end{itemize} }

The proof will establish the lemma for any Moebius transformation
defined on $T(z)$. We have associated in the last section to
every $(n-2)$-sphere $S_z$ an euclidean center $x$ and a radius vector
 $y$ where $z=x+\i y$. The lemma asserts that $
\gamma(S_z)=S_{\tilde \gamma(z)}$, if $\gamma$ is orientation
preserving. That result could be a good motivation for the
introduction of complex points in elementary
inversive geometry.

\dem The lemma is obvious for euclidean isometries and
homotheties. Let's prove it for any inversion $\gamma$. We take for $\gamma$ the inversion through the unit
sphere of $\R^n$,  We have $\gamma (x)={x\over \|x\|^2}={x\over
\Q(x)}$. If $\Q(z)\Q(z')\not=0$,
$$\Q(\tilde\gamma (z)-\tilde\gamma(z'))=\Q\left[{z\Q(z')-z'\Q(z)\over \Q(z)\Q(z')}\right]={\Q(z-z')\over \Q(z)\Q(z')}$$
Therefore, $\Q(z-z')=0 \Leftrightarrow \Q\left(\tilde\gamma(z)-\tilde\gamma(z')\right)=0$. If we take $z'$ real, and write $z=x+\imath y$, the last equation
says that setwise,
  $\gamma(T(z))$ and $T(\tilde \gamma(z))$ are equal, which means $\gamma(T(z))
=T(\tilde \gamma(x+\imath \varepsilon(z,\gamma)y))$, where $\varepsilon (z,\gamma)=\pm
1$. As previously, by the same topological arguments, one easily sees that $\varepsilon =1$ for orientation preserving maps and $-1$ for the others.
The lemma and the proposition are proved.

\Th{Theorem}{The fourth symmetric domain of Elie Cartan, called the Lie ball, is equivariantly  diffeomorphic to $T\Bn$,
and the action of the hyperbolic group extends to its boundary.
So, the closure $\overline{\B^n}$ of  Lie ball is an
equivariant compactification of the tangent bundle to 
hyperbolic space.}

\dem The equivariance is obvious. We prove that
$\theta=T^{-1}\circ S$ is a diffeomorphism, by direct computation.
This allows us to pull back on the tangent space $T\Bn$ the Bergmann metric
of $\B^n$. We use the chart
$\Bn\times \R^n$ on $T\Bn$.

 Let $v$ be a tangent vector at point $x$. To
avoid confusion, we will sometimes denote $v$ by $(x,v)$. we have
$ \theta(v)=z$ such that $T(z)=S(v)$. We recall that the euclidean and  the hyperbolic centers of spheres
coincide at the origin. The isometry $\delta_x$ is an
 involution sending $x$ to the origin, and $v$ to $v'$, by the tangent map. We have
\begin{eqnarray*}
  \theta (v) &=& T^{-1}\circ S\circ (d\delta_x(o)\circ d\delta_x(x))(v) \\
   &=& T^{-1}\circ  \delta_x\circ S\circ d\delta_x(x)(v) \\
  &=& \tilde\delta_x\circ (T^{-1}\circ S)(d\delta_x(x)(v))\\
  &=& \tilde\delta_x\circ\theta(v'), \quad \mbox{where $v'=d\delta_x(x)(v)={v\over \|x\|^2-1}\cdot$}
\end{eqnarray*}
The euclidean
radius of $S(v')$ (centered at the origin) is $\th\|v'\|$,\\
 so \quad $\theta (v')=\theta (o,v')=\i{v'\over
\|v'\|}\th\|v'\|=-\i{v\over\| v\|}\th{\|v\|\over 1- \|x\|^2}$.\\
The maps in this
construction are differentiable of maximum rank, as can be seen
for $\theta (v')$ by limited expansion at order 1 in a
neighborhood of $v=0$. The final expression of $\theta$ is:

\begin{eqnarray}
  \theta (x,v)=\tilde\delta_x(z')&=& {(\|x\|^2-1)z'+(1+\Q(z'))x-2(z'\cdot x)x\over \Q(x)\Q(z')-2x\cdot
  z'+1},\\
   \mbox {with}\quad z'=\theta (o,v')&=& -\i{v\over\| v\|}\th{\|v\|\over 1-
   \|x\|^2}\nonumber
\end{eqnarray}

The real part and the imaginary part of $\theta (x,v)$ are respectively the
euclidean center and the radius vector of the image $(n-2)$-sphere.

As an exercise, we derive the
following well-known result:
\Th{corollary}{Let $S(c,r)$ be a
sphere in the Poincaré model of the hyperbolic space, with
hyperbolic center $c$ and hyperbolic radius $r$. Its euclidean
center $c_e$ and radius $r_e$ are given by
$$c_e={(1-\alpha^2)c   \over 1- \|c\|^2\alpha^2},\quad
r_e={(1-\|c\|^2)\alpha \over 1- \|c\|^2\alpha^2}, \
\mbox{where} \ \ \alpha= \th{r\over 1-\|c\|^2}$$}
 \dem Let's
consider the hyperbolic space $\mathbb{B}^n$ as a coordinate hyperplane
of $\mathbb{B}^{n+1}$. Its spheres can be seen as
$(n-2)$-spheres of $\mathbb{B}^{n+1}$. The spheres
$S(c,r)$ are represented by vectors orthogonal to $\Bn$.
The calculation follows.

\subsubsection*{The symmetric Kähler structure on the tangent bundle}
In the
coordinate system $(x_1, \ldots, x_n, x_{n+1}, \ldots,x_{2n})$ with
$x_{n+i}={\partial \over \partial x_i}$, we carry by $\theta$ the metric of the Lie ball on $T\Bn$. Let's first compute the
metric on the tangent space at the origin, and then carry
it everywhere with the tangent maps to isometries $\delta_x$ of the base
manifold.

On the tangent bundle $TT\Bn$ to $T\Bn$, we use the coordinate system:\\
 $(y_1,\ldots, y_{2n},{\partial \over
\partial y_1}, \ldots,{\partial \over \partial y_{2n}})$. Let
$v\equiv (0,v)$ be a tangent vector at the origin in $\Bn$. Let
$(v,w)\in TT\Bn$ be a tangent vector at $v$, with $w=(u,u')$, where
$u=\sum_{1}^{n}w_i {\partial \over \partial y_i}$ and
$u'=\sum_{n+1}^{2n}w_i {\partial \over \partial y_i}$. The
differential $d\theta$ near $x=0$ can be obtained by
limited expansion of expression (1) at order one. We obtain:
$$\theta (x,v)= x-\left(x-2\left(x\cdot {v\over \|v\|}\right){v\over
\|v\|}\right)\th^2\|v\|+\imath {v\over \|v\|}\th\|v\|+o(x).$$
 That
expression has a principal part whose derivative at $(0,v)$
preserves the real and imaginary subspaces, with a Jacobian matrix
$d\theta (v)=\left(\begin{array}{cc}A_1&0\\0&A_2\end{array}\right)$,
where $A_1$ and $A_2$ are matrix of order $n$.
The Bergmann metric at $\theta(o,v)$ in the Lie ball is isotropic: in fact, for $x=o$, $\theta (v)=\i y$ where $y= {v\over \|v\|}\th\|v\|$
is real. Let  $B$ be any orthonormal basis of $T\Bn$ at the origin
$o$. Since the metric at $o$ is preserved by the action of its
stabilizer $O(n)\times SO(2)$, all the
vectors of the basis $(B, \imath B)$ have the same size, which
proves that the Bergmann metric is isotropic at $o$, then equal to
$\sum_{i}^{n}\|dz_i\|^2$, up to a positive factor. The holomorphic
continuation of the hyperbolic motion $\delta_y$ sends $o$ to $y$,
so that the metric at $y$ is ${1\over (1-
\|y\|^2)^2}\sum_{1}^{n}\|dz_i\|^2={\cosh^4\|v\|}\sum_{1}^{n}\|dz_i\|^2$.
The Bergmann metric keeps the same value at $\i y$. Finally, the
metric at $(o,v)\in T\Bn$ is
$$ds^2(w)=
{^tu}B_1u + {^tu'}B_2u', \quad\mbox{where}\
B_i=\cosh^4\|v\|\,{^t{\!\!A}}_iA_i, \ i=1,2.$$

It follows that at $(o,v)$ in $TT\Bn$, the vector subspace $u=0$
tangent to the fiber is orthogonal with respect to this new
riemannian metric to the subspace $u'=0$ which can be understood
as tangent to a 'lift' of the base manifold $\Bn$. The theory of connections
makes this notion more precise: given a riemannian metric $g$ on a
manifold, one can deduce a unique torsion free associated affine
connexion called the Levi-Civita connexion.

 Let $V$ be a manifold, and
$\nabla$ a connection on $V$. In the second tangent space $TTV$,
one defines the vertical space at a point $(x,v)\in TV$ as the
subspace tangent to the fiber. Let $H$ be the set of tangent
vectors to $V$ parallel to $(x,v)$, in a neighborhood of $(x,v)$
where parallel transportation is well defined. $H$ is a manifold
in a neighborhood of $(x,v)$,
 The horizontal space is the set of vectors tangent to $H$ at $(x,v)$.
That decomposition is preserved by any isometry of the connexion $\nabla$.

In our case, the two spaces are orthogonal
at $(o,v)$, so they are orthogonal at every tangent point, since
the isometry group is transitive on $\Bn$. In such a case,
the geodesic flow is incompressible [1].

This new metric has interesting properties. 

\Th{Theorem}{The metric on the tangent space $TM$ has the following properties:
\begin{itemize}
	\item Every vector line is a geodesic.
	\item The tangent space at a point of $M$ is totally geodesic.
	\item The tangent space to a geodesic of $M$, also called a leaf of the riemannian foliation, is a hyperbolic disk.
\end{itemize}}
\dem In the Lie ball, every straight
line through the origin is a geodesic. Thus, every vector line in the
tangent space is therefore a geodesic.

The tangent space at every point $p$ is totally geodesic: we suppose that $p$ is the origin. The image of the tangent space to $p$ by $\theta$ is $\imath \Bn$, which is the image of $\Bn$ by an isometry of $\B^n$.
But $\Bn$ is totally geodesic.

The tangent space to a geodesic is of special interest. We suppose
that the geodesic goes through the origin. It is a straight line $l$
in the Poincaré model, directed by a unit vector $\vec{u}$. A tangent vector to a point of $l$ is
represented by a $(n-2)$-sphere centered on $l$ at point $\alpha$,
with radius $\beta$, contained in the unit sphere. So its complex
coordinate in the Lie ball is $\alpha+\imath \beta \vec{u}$ where
$\|\alpha\|^2+\|\beta\|^2<1$.

Then, the tangent space to a geodesic is a hyperbolic disk. It is
not hard to show that its stabilizer is a normal extension of the orientation preserving hyperbolic group of the disk.

Every complete hyperbolic manifold $M$ is a quotient manifold of $\mathbb{B}^n$ by a discrete group
of Moebius transformations. Its tangent space is therefore a quotient manifold of the
Lie ball. We thus define a compactification of $TM$.

 A deep study of the geometry of the Lie ball seems useful. It will give a better
 insight of
 the hyperbolic space and help to derive new results through the compactification of tangent spaces to hyperbolic manifolds by  codimension 2 horospheres.

The above study is being  adapted for a few other
riemannian symmetric spaces of non compact type.\bigskip

\parindent -5pt
{\centering \large \bf REFERENCES\par }\penalty 10000
\medskip
\small

[1] M. T. K. Abbassi, M. Saari, On natural metrics on tangent bundles of riemannian manifolds, {\it Archivum Mathematicum}, Brno,
Tomus \textbf{41} (2005),71-92

[2] D. N. Akhiezer,  S. G. Gindikin, On Stein extension of real symmetric spaces, {\it Math Annalen}, \textbf{286}, (1990) n° 1-3,
1-12

[3] V. Avanissian, {\it Cellule d'harmonicit\'e et prolongement analytique complexe}, {Coll. Travaux en cours}
Hermann, 1985.

[4] M. Berger, {\it Espaces euclidiens, triangles, cercles et sph\`eres}, Cedic/Nathan, 1977

[5] S. Bergmann,  \"Uber Die Kernfunction eines Bereiches und ihr Verhalten am Rande, {\it J. Reine Angew. Math.,}
\textbf{169}, 1--42.

[6] M. Boutaleb \emph{Sur la cellule d'harmonicit\'{e} de la boule
unit\'{e} de $\mathbb{R}^n$}, Doctorat de 3$^e$cycle, U.L.P. Strasbourg, France (1983).

[7] D. Burns, S., Halverscheid, R., Hind,  The geometry of Grauert tubes and complexification of
symmetric spaces, \emph{Arxiv math} CV/0109186v1, sept 8, 2005

[8] E. Cartan, Sur les domaines born\'es homog\`enes de l'espace de $n$ variables complexes, {\it Abh. Math. Sem.
Ham. Univ.} \textbf{11} (1935) 116-162.

[9] V. Guillemin, M. Stenzel,  Grauert tubes and the homogenous Monge-Ampère equation, \emph{J. Diff. Geom.},
\textbf{34} (1991), 561-570

[10] Kol\'{a}r, J., P. W. Michor, Slov\'{a}k, \emph{J., Natural operations in differential geometry},
Springer-Verlag, Berlin, 1993

[11] O. Kowalski, M. Sekizawa, Natural transformations of riemannian metrics on manifolds to metrics on tangent bundles
- a classification, \emph{Bull. Tokyo Gakugei Univ.} (4) 40 (1988) 1-29

[12] P. Lelong, Prolongement analytique et singularit\'es complexes des fonctions harmoniques,
{\emph{Bull. Soc. Math. Belg.}} \textbf{7} (1954-55), 10--23.

[13] L. Lempert, R. Szöke, Global solutions of the homogeneous complex Monge-Ampère equation
and complex structures on the tangent bundle of Riemannian manifolds, \emph{Math. Annalen}, \textbf{290} (1991), 689-712

[14] V. Oproiu, A locally symmetric Kähler Einstein structure on the tangent bundle of a space form,
\emph{Beiträge zur Algebra und Geometrie / Contributions to Algebra and Geometry} \textbf{40} (1999) 363-372

[15] W. A. Thurston, \emph{Three-Dimensional Geometry and Topology}, vol. \textbf{1}, Princeton University Press, Princeton, N.J. (1997)

[16] J. A. Wolf, R. Zierau, Linear cycle spaces in flag domains, \emph{Math. Annalen},  \textbf{316} (2000), 529-545.

\vskip 2mm
\scriptsize
\parindent -8pt
\begin{tabular}{l@{\qquad\qquad}l}
\it Address : \rm Tchangang Tambekou Roger&\it  Teaching address :\\
\verb"tchangang_tambekou@yahoo.com"&Département de Mathématiques \\
Centre de Recherches en Géométrie et Applications&Faculté des Sciences\\
B.P. 8 451 Yaoundé&Université de Yaoundé I\\
Cameroun&Cameroun
\end{tabular}

\end{document}